\documentclass[12pt]{amsart}
\usepackage{amscd,amssymb}
\usepackage[arrow,matrix]{xy}
\usepackage[plainpages,backref,urlcolor=blue]{hyperref}

\theoremstyle{plain}
\newtheorem{thm}[subsection]{Theorem}

\newtheorem{prop}[subsection]{Proposition}
\newtheorem{cor}[subsection]{Corollary}

\theoremstyle{definition}
\newtheorem{rk}[subsection]{Remark}

\newtheorem{ex}[subsection]{Example}

\numberwithin{equation}{section}
\newcommand{\OO}{{\mathcal O}}

\newcommand{\A}{{\mathcal A}}

\newcommand{\al}{{\alpha}}
\newcommand{\be}{{\beta}}
\newcommand{\lb}{{\lambda}}
\newcommand{\ga}{{\gamma}}

\newcommand{\Z}{\mathbb{Z}}
\newcommand{\Q}{\mathbb{Q}}

\newcommand{\C}{\mathbb{C}}

\newcommand{\PP}{\mathbb{P}}

\begin{document}

\title [Poincar\'e-Deligne polynomial of triple point line arrangements]
{The Poincar\'e-Deligne polynomial of Milnor fibers of triple point line arrangements is combinatorially determined}

\author[Alexandru Dimca]{Alexandru Dimca$^1$}
\address{Univ. Nice Sophia Antipolis, CNRS,  LJAD, UMR 7351, 06100 Nice, France.}
\email{dimca@unice.fr}

\thanks{$^1$ Partially supported by Institut Universitaire de France} 



\begin{abstract} Using a recent result by S. Papadima and A. Suciu, we show that 
the equivariant Poincar\'e-Deligne polynomial  of the Milnor fiber of a projective line arrangement having only double and triple points is combinatorially determined.

\end{abstract}

\maketitle


\section{Introduction} \label{sec1}

Let $\A$ be an arrangement of $d$ hyperplanes in $\PP^{n}$, with $d \geq 2$, given by a reduced equation 
$ Q(x)=0$. 
Consider the corresponding complement $M$ defined by $Q(x)\ne 0$ in $\PP^{n}$, and the global Milnor fiber $F$ defined by $Q(x)-1=0$ in $\C^{n+1}$ with monodromy action $h:F \to F$, $h(x)=\exp(2\pi i/d)\cdot x$. We refer to \cite{OT} for the general theory of hyperplane arrangements.

The following are basic open questions in this area, a positive answer for any question in this list implying the same for the previous ones.

\begin{enumerate}
\item Are the Betti numbers $b_j(F)$  combinatorially determined, i.e. determined by the intersection lattice $L(\A)$ ?
\item Are the monodromy operators $h^j:H^j(F) \to H^j(F)$  combinatorially determined?
\item Is the equivariant Poincar\'e-Deligne polynomial $PD^{\mu_d}(F)$ of $F$ coming from the monodromy action combinatorially determined? Here $\mu_d$ is the multiplicative group of $d$-th roots of unity and the definition of $PD^{\mu_d}(F)$ is recalled in the next section.
\end{enumerate}

On the positive side, it was shown by N. Budur and M. Saito in \cite{BS} that the spectrum $Sp(\A)$ of $\A$, whose definition is also recalled in the next section, is combinatorially determined.

\medskip

We assume in the sequel that $n=2$ and that the line arrangement $\A$ has only double and triple points. Then a recent result of S. Papadima and A. Suciu \cite{PStrip} shows that the answer to the question (2) above is positive. More precisely, they have introduced a combinatorial invariant $\beta_3(\A)\in \{0,1,2\}$ of the line arrangement $\A$ such that the multiplicity of a cubic root of unity $\lambda \ne 1$ as an eigenvalue for $h^1$ is exactly $\beta_3(\A)$.

The main result of this note, answering a question raised by Alex Suciu, is the following.

\begin{thm} \label{Thm1} 
Let $\A$ be an arrangement of $d$ lines in $\PP^{2}$,  such that $\A$ has only double and triple points. Then the equivariant Poincar\'e-Deligne polynomial $PD^{\mu_d}(F;u,v,t)$ of $F$ coming from the monodromy action is determined by  the number $d$ of lines in $\A$, the number $n_3(\A)$ of triple points in $\A$  and the Papadima-Suciu invariant $\beta_3(\A)$.

\end{thm}
In particular, the question (3) above has a positive answer in this case. This is rather surprising, given the complexity of the mixed Hodge structure on the cohomology of the Milnor fiber $F$, as seen from Propositions \ref{prop1} and \ref{prop2} below. The very explicit formulas given in these two Propositions apply to certain monodromy eigenvalues for arbitrary line arrangements, see Remarks \ref{rk3} and \ref{rk4} below.

For a possible application to the study of some (singular) projective surfaces, see Remark \ref{rk1}. Relations to the superabundance or the defect of some linear systems passing through the triple points of $\A$ are described in Remark \ref{rk2}.

Note also that there are very few examples of nonisolated (quasi-homogeneous)  hypersurface singularities $(X,0)$ for which both the monodromy and the MHS on the corresponding Milnor fibers are well understood, even though the isolated quasi-homogeneous case was settled by  J. Steenbrink \cite{St} a long time ago. 

The case of a hyperplane arrangement in $\PP^{3k-1}$ which is obtained by taking a product $\A_1 \times \A_2 \times ...\times \A_k$ of $k$ line arrangements $\A_j$ 
having  only double and triple points can now be treated using the results in this note and Theorem 1.4 in \cite{DNag}.

In the last section we prove the following related result.

\begin{prop} \label{prop3} 
Let $C:Q=0$ be a degree $d$ reduced curve in the projective plane $\PP^{2}$,  and let $F:Q-1=0$ be the associated Milnor fiber in $\C^3$. Then the equivariant Poincar\'e-Deligne polynomial $PD^{\mu_d}(F;u,v,t)$ of $F$ coming from the monodromy action is determined by its specialization, the Hodge-Deligne polynomial
$$HD^{\mu_d}(F;u,v)=PD^{\mu_d}(F;u,v,-1).$$
\end{prop}
Since the Hodge-Deligne polynomial (or rather a compactly supported version of it, is additive, see for instance \cite{DL12}), this result might be used in some situations to compute these polynomials.
It is an {\it open question} whether such a result holds in higher dimensions, even for the hyperplane arrangements. 

For similar non-cancellation properties in the case of braid arrangements $A_3$ and $A_4$, see \cite{DL13}, section 6.

\section{Equivariant Hodge-Deligne and Poincar\'e-Deligne  polynomials and spectra}

Recall that if $X$ is a quasi-projective variety over $\C$ one can consider the Deligne mixed Hodge structure (for short MHS) on the rational cohomology groups
$H^*(X,\Q)$ of $X$. We refer to \cite{PS} 
for general notions and results concerning the MHS.

Since this MHS is functorial with respect to algebraic mappings, if a finite group $\Gamma$ acts
algebraically on $X$, each of the graded pieces
\begin{equation} 
\label{eq:hpq}
H^{p,q}(H^j(X,\C)):=Gr_F^pGr^W_{p+q}H^j(X,\C)
\end{equation}
becomes a $\Gamma$-module in the usual functorial way, and these modules are the building 
blocks of the Hodge-Deligne polynomial $HD^{\Gamma}(X;u,v)\in R(\Gamma)[u,v]$, defined by
\begin{equation} 
\label{eq:hd}
HD^{\Gamma}(X;u,v)=\sum_{p,q}E^{\Gamma;p,q}(X)u^pv^q,
\end{equation}
where $E^{\Gamma;p,q}(X)=\sum_j(-1)^jH^{p,q}(H^j(X,\C))\in R(\Gamma)$. 
One may consider an even finer (and hence harder to determine) invariant, namely the equivariant Poincar\'e-Deligne polynomial
\begin{equation} 
\label{eq:pd}
PD^{\Gamma}(X;u,v,t)=\sum_{p,q,j} H^{p,q}(H^j(X,\C))u^pv^qt^j \in R_+(\Gamma)[u,v,t].
\end{equation}
Clearly one has $PD^{\Gamma}(X;u,v,-1)=HD^{\Gamma}(X;u,v)$.
 
The case of interest to us is when $\Gamma=\mu_d$ and the action on $F$ is determined by
$$\exp(2\pi i/d) \cdot x=h^{-1}(x).$$
The reason to use $h^{-1}$ instead of $h$ is either functorial (i.e. to really have a group action when $\Gamma$ is not commutative, see \cite{DL13}) or geometrical, as explained in \cite{DS14}, in order to get results compatible with those in \cite{BS}, which we survey below.
Recall that the spectrum of a hyperplane arrangement $\A \subset \PP^n$ is the polynomial $$Sp(\A)= \displaystyle{\sum_{\alpha \in \Q} n_{\alpha}\,\, t^{\alpha}},\,$$
with coefficients given by
$$n_{\alpha} = \displaystyle{\sum_j (-1)^{j-n} \dim\,Gr_F^p \tilde{H}^j(F,\C)_{\lambda}},\,$$
where $\,\,\,p= \lfloor n +1- \alpha \rfloor,\,\,\lambda= exp(-2i\pi \alpha),$
with $\tilde{H}^j(F,\C)_{\lambda} = H^j(F,\C)_{\lambda}$ (eigenspaces with respect to the action of $(h^j)^{-1}$ as explained above) for $ j \neq 0$,   
 $\tilde{H}^0(F,\C)_{\lambda}=0$
and $\lfloor y \rfloor := max \{k\in \Z\,|\, k \leq y \}.$ It is easy to see that $n_{\alpha}=0$ for $\alpha \notin (0,n+1)$.

Theorem 3 in \cite{BS} implies the following result.

\begin{thm} \label{jumpcoef}
Let $\A$ be an arrangement of $d$ lines in $\PP^{2}$,  such that $\A$ has only double and triple points. Let $n_3(\A)$ denote the number of triple points in $\A$. Then 
 $n_{\alpha} = 0$ if $\alpha d \notin \Z,\,$ and for $\alpha = \frac{j}{d} \in\, ]0,1]\,$ with $j \in [1,d] \cap \Z$ the following hold.
\begin{center}
$n_{\alpha} = \binom {j-1} {2} - \displaystyle{n_3(\A) \binom{ \lceil 3j /d \rceil -1} {2}} $\\
$n_{\alpha +1} = (j-1)(d-j-1) - \displaystyle{n_3(\A)(\lceil 3j /d \rceil -1)(3- \lceil 3j /d \rceil )}$ \\
$n_{\alpha+2} = \binom {d-j-1} {2} - \displaystyle{ n_3(\A)\binom{ 3-\lceil 3j /d \rceil } {2}} - \delta_{j,d}$
\end{center}
where $\lceil y \rceil:= min\{k\in \Z\,|\, k\geq y\},\,$ and $\delta_{j,d}=1$ if $j=d$ and $0$ otherwise.

In particular, the spectrum $Sp(\A)$ is determined by the number $d$ of lines in $\A$ and the number $n_3(\A)$ of triple points. 
\end{thm}

\section{The proof of Theorem \ref{Thm1} }
Let us consider the cohomology groups $H^j(F,\Q)$ one by one and discuss the corresponding MHS and monodromy action. The group $H^0(F,\C)$ is clearly one dimensional, of type $(0,0)$ and the monodromy action is trivial, i.e. $H^0(F,\C)=H^0(F,\C)_1$.

The next group $H^1(F,\Q)$ is already more interesting. It has a direct sum decomposition 
$$H^1(F,\Q)=H^1(F,\Q)_1 \oplus H^1(F,\Q)_{\ne 1}$$
in the category of MHS. The fixed part under the monodromy, $H^1(F,\Q)_1$ is isomorphic to the cohomology group of the projective complement, namely  $H^1(M,\Q)$ and hence it has dimension $d-1$ and it is a pure Hodge-Tate structure of type $(1,1)$.

The other summand $H^1(F,\Q)_{\ne 1}$ is always a pure Hodge structure of weight 1, see
 \cite{BDS} and \cite{DP}
for two distinct proofs. Moreover, in the case when only double and triple points occur in $\A$, then the second summand corresponds to cubic roots of unity and it can be non zero only when $d$ is divisible by 3, see for instance Remark \ref{rk3} below. By combining Papadima-Suciu results in \cite{PStrip}  with  (the proof) of Theorem 1 in \cite{Dtrip} (see also Theorem 2 in \cite{BDS} for a more general result and Remark \ref{rk2} below for additional information), one gets
\begin{equation} 
 \label{H1}
h^{1,0}(H^1(F))_{\ga'}=h^{0,1}(H^1(F))_{\ga}=\beta_3(\A) \text{ and } h^{1,0}(H^1(F))_{\ga}=h^{0,1}(H^1(F))_{\ga'}=0,
\end{equation}
where $\be=1/3$, $\ga =\exp(-2\pi i\be)$, $\be'=2/3$, $\ga' =\exp(-2\pi i \be')=\overline \ga$. Here and in the sequel we use the notation
$h^{p,q}(H^j(F))_{\lb}$ to denote the multiplicity of the 1-dimensional $\mu_d$-representation $r_{\lb}$  sending $\exp(2\pi i/d)$ to $\lb \in \mu_d \subset \C^*=Aut(\C)$ in the $\mu_d$-module
$H^{p,q}(H^j(F,\C))$ defined in \eqref{eq:hpq}. Note that one has

\begin{equation} 
\label{relhpq}
\dim\,Gr_F^p {H}^j(F,\C)_{\lambda}=\sum_{q\geq j-p}h^{p,q}(H^j(F))_{\lb},
\end{equation}
by the general properties of MHS, $F$ being smooth.

To determine the  equivariant Poincar\'e-Deligne polynomial $PD^{\mu_d}(F)$ of $F$ is clearly equivalent to determine all the equivariant mixed Hodge numbers $h^{p,q}(H^j(F))_{\lb}$. Until now, we have done this for $j=0$ and $j=1$.

It remains to treat the case $j=2$, which is the most difficult.  We have again a direct sum decomposition 
$$H^2(F,\Q)=H^2(F,\Q)_1 \oplus H^2(F,\Q)_{\ne 1}$$
in the category of MHS. The fixed part under the monodromy, $H^2(F,\Q)_1$ is isomorphic to the cohomology group of the projective complement, namely  $H^2(M,\Q)$ and hence has dimension $b_2(M)$ and pure Hodge-Tate type $(2,2)$. Since the Euler characteristic
$\chi(M) =b_0(M)-b_1(M)+b_2(M)$ can be computed from the combinatorics, it follows that
\begin{equation} 
\label{H2_1}
b_2(M)=\displaystyle{ \binom{ d -1} {2}} -n_3(\A).
\end{equation}
We can also write $H^2(F,\Q)_{\ne 1}$ as a direct sum of two MHS, namely $H^2(F,\Q)_{\ne 1}=H \oplus H'$, where $H$ corresponds to the eigenvalues of $h^2$ which are cubic roots of unity different from 1, and $H'$ corresponds to all the other eigenvalues.

Proposition 4.1 in \cite{DNag} implies that $H'$ is a pure Hodge structure of weight 2, i.e.
$h^{p,q}(H^2(F))_{\lb}=0$ for $p+q \ne 2$ and $\lb$ not a cubic root of unity.
On the other hand, Theorem 1.3 in \cite{DL12} implies that the only weights possible for $H$ are 2 and 3, hence $h^{p,q}(H^2(F))_{\lb}=0$ for $p+q \notin \{ 2,3\}$ and $\lb$  a cubic root of unity.

Now we explicitly determine the equivariant mixed Hodge numbers $h^{p,q}(H^2(F))_{\lb}$
for $\lb \ne 1$, the case $\lb = 1$ being already clear by the above discussion.
The above discussion implies also the following result.

\begin{prop} 
\label{prop1} Let $\A$ be an arrangement of $d$ lines in $\PP^{2}$,  such that $\A$ has only double and triple points. Let $n_3(\A)$ denote the number of triple points in $\A$.
Assume that $\lb=\exp(-2 \pi \al)$, with $0<\al=j/d<1$, is not a cubic root of unity. Then one has $h^{2,0}(H^2(F))_{\lb}=n_{\al}$, $h^{1,1}(H^2(F))_{\lb}=n_{\al+1}$ and $h^{0,2}(H^2(F))_{\lb}=n_{\al+2}$, where  the integers $n_{\al}, n_{\al+1}, n_{\al+2}$ are given by the formulas in Theorem \ref{jumpcoef}.

\end{prop} 

\begin{rk} 
\label{rk3} Let $\A$ be any essential arrangement of $d$ lines in $\PP^{2}$, i.e. $\A$ is not a pencil of lines through a point. Then the formulas given in Proposition \ref{prop1} hold for any $\lb \in \mu_d$ such that there is a line $L \in \A$ with $\lb ^k \ne 1$ whenever there is a point of multiplicity $k$ in $\A$ situated on $L$. Indeed, this last condition is known to imply that $H^1(F)_{\lb}=0$, see \cite{Li2}. In such a case,  the integers $n_{\al} $ are not given by the formulas in Theorem \ref{jumpcoef}, but they are described in Theorem 3 in \cite{BS}.

\end{rk}

Now we consider the case of the cubic roots of unity $\ga=\exp(-2\pi i\be)$ and $\ga'=\exp(-2\pi i\be')$ introduced above. They can be eigenvalues of $h^2$ only when $d$ is divisible by 3.

\begin{prop} 
\label{prop2} Let $\A$ be an arrangement of $d$ lines in $\PP^{2}$,  such that $\A$ has only double and triple points. Let $n_3(\A)$ denote the number of triple points in $\A$ and suppose that $d$ is divisible by 3.
Then one has the following.

\begin{enumerate}
\item $h^{2,0}(H^2(F))_{\ga}=h^{0,2}(H^2(F))_{\ga'}  =n_{\be'+2}$;

\item $h^{1,1}(H^2(F))_{\ga}=h^{1,1}(H^2(F))_{\ga'}  =n_{\be'+2}+n_{\be'+1} -n_{\be} +\be_3(\A)$;

\item $h^{0,2}(H^2(F))_{\ga}=h^{2,0}(H^2(F))_{\ga'}  =n_{\be'+2}+n_{\be'+1} +n_{\be'}-n_{\be}-n_{\be+1} +\be_3(\A)$;

\item $h^{2,1}(H^2(F))_{\ga}=h^{1,2}(H^2(F))_{\ga'}  =n_{\be}-n_{\be'+2}$;

\item $h^{1,2}(H^2(F))_{\ga}=h^{2,1}(H^2(F))_{\ga'}  =n_{\be+1}+n_{\be}-n_{\be'+1}-n_{\be'+2}-\be_3(\A)$.

\end{enumerate}
Here $\be=1/3$, $\be'=2/3$ and the various integers $n_{\eta}$ are given by the formulas in Theorem \ref{jumpcoef}.

\end{prop} 

\proof

In the case $\al=\be$, the definition of the spectrum and the above discussion on the mixed Hodge properties of the cohomology group of the Milnor fiber $F$ yield the following relations.

\begin{enumerate}
\item $n_{\be}= h^{2,0}(H^2(F))_{\ga}+h^{2,1}(H^2(F))_{\ga}$;

\item $n_{\be+1}= h^{1,1}(H^2(F))_{\ga}+h^{1,2}(H^2(F))_{\ga}$ (use \eqref{H1});

\item $n_{\be+2}= h^{0,2}(H^2(F))_{\ga}-h^{0,1}(H^1(F))_{\ga}=h^{0,2}(H^2(F))_{\ga}-\be_3(\A)  $ (use again \eqref{H1});

\end{enumerate}
Similarly, for $\al=\be'$, we get the following.

\begin{enumerate}
\item $n_{\be'}= h^{2,0}(H^2(F))_{\ga'}+h^{2,1}(H^2(F))_{\ga'}$;

\item $n_{\be'+1}= h^{1,1}(H^2(F))_{\ga'}+h^{1,2}(H^2(F))_{\ga'}-\be_3(\A) $ (use \eqref{H1});

\item $n_{\be'+2}= h^{0,2}(H^2(F))_{\ga'}$ (use again \eqref{H1});

\end{enumerate}
The proof is completed by using the obvious equality 
$$h^{p,q}(H^2(F))_{\lb}=h^{q,p}(H^2(F))_{\overline \lb},$$
obtained by taking the complex conjugation.

\endproof

\begin{rk} 
\label{rk4} Let $\A$ be any essential arrangement of $d$ lines in $\PP^{2}$, i.e. $\A$ is not a pencil of lines through a point. Then the formulas given in Proposition \ref{prop2} where we take $\be_3(\A)=0$ clearly  hold for any $\lb \in \mu_d$ such that $H^1(F)_{\lb}=0$, with the integers $n_{\al} $  given by  Theorem 3 in \cite{BS}.

\end{rk} 

Moreover, it is clear that Propositions \ref{prop1} and \ref{prop2} imply our Theorem \ref{Thm1}.
They also yield the following.

\begin{cor} 
\label{cor1} Let $\A$ be an arrangement of $d$ lines in $\PP^{2}$,  such that $\A$ has only double and triple points. Then the dimensions $\dim Gr^W_2H^2(F,\Q)$ and $\dim Gr^W_3H^2(F,\Q)$ of the graded quotients with respect to the weight filtration $W$ depend both on the Papadima-Suciu invariant $\be_3(\A)$.

\end{cor}

\begin{ex} 
\label{ex1} Let $\A$ be the Ceva (or Fermat) arrangement of $d=9$ lines in $\PP^{2}$ given by the equation
$$Q(x,y,z)=(x^3-y^3)(x^3-z^3)(y^3-z^3).$$ 
Then the Papadima-Suciu invariant $\be_3(\A)$ is equal to 2, there are $n_3(\A)=12$ triple points and a direct computation gives the following formula for the spectrum
$$Sp(\A)= t^{1/3} +3 t^{4/9} +6t^{5/9} +10t^{2/3} +3t^{7/9} +9t^{8/9} +16t +$$
$$+6t^{11/9} +10t^{4/3} -2t^{5/3} +6t^{16/9} -8 t^2 
+9t^{19/9} +3t^{20/9} -2t^{7/3} +6t^{22/9} +3t^{23/9} +t^{8/3} -t^3 .$$
Proposition \ref{prop2} implies
$$h^{2,1}(H^2(F))_{\ga}=h^{1,2}(H^2(F))_{\ga'}  =n_{1/3}-n_{8/3}= 1-1=0$$
and
$$h^{1,2}(H^2(F))_{\ga}=h^{2,1}(H^2(F))_{\ga'}  =n_{4/3}+n_{1/3}-n_{5/3}-n_{8/3}-\be_3(\A)=10+1+2-1-2=10.$$
These values correct a misprint in \cite{DL12}, p. 244 and confirm the computations done by P. Bailet in \cite{Bai}. This examples also shows that the integers $n_{\eta}$ may be strictly negative.
\end{ex} 

\begin{rk} 
\label{rk1} Let $\A$ be an arrangement of $d$ lines in $\PP^{2}$,  such that $\A$ has only double and triple points. Then, in view of Theorem 1.1 in \cite{DL12}, the results in Propositions \ref{prop1} and \ref{prop2} can be used to describe the $\mu_d$-action on the cohomology of the associated surface 
$$X_Q: Q(x,y,z)-t^d=0$$
in $\PP^3$, which is a singular compactification of the Milnor fiber $F$.

\end{rk}

\begin{rk} 
\label{rk2} Let $\A$ be an arrangement of $d$ lines in $\PP^{2}$,  such that $\A$ has only double and triple points and $d=3m$ for some integer $m$. Let $T \subset \PP^2$ be the set of triple points in $\A$. If $S=\C[x,y,z]$ denotes the graded ring of polynomials in $x,y,z$, consider  the evaluation map
\begin{equation} \label{eval}
\rho: S_{2m-3} \to \C^T
\end{equation}
obtained by picking up a representative $s_t$ in $\C^3$ for each point $t\in T$ and sending
a homogeneous polynomial $h \in S_{2m-3}$ to the family $(h(s_t))_{t\in T}$.

Then Theorem 2 in \cite{BDS} and the discussion folowing it imply the key formula \eqref{H1} above. This can be reformulated as
\begin{equation} \label{eval2}
\be_3(\A)=\dim (Coker \rho),
\end{equation}
and the last integer is by definition  the {\it superabundance}  or the {\it defect} $S_{2m-3}(T)$ of the finite set of points $T$ with respect to the polynomials in $S_{2m-3}$. Since by the work of S. Papadima and A. Suciu we know that $\be_3(\A) \in \{0,1,2\}$, this gives a very strong restriction on the position of the triple points in such a line arrangement. For other relations to algebraic geometry of a similar flavor we refer to \cite{E}, \cite{LV}, \cite{Li}.

\end{rk} 

\section{The proof of Proposition \ref{prop3} }
We follow the notation from the previous section, in particular $M$ denotes the complement of $C$ in $\PP^2$ given by $Q \ne 0$.
To prove Proposition \ref{prop3}, we have to check whether for each character $r_{\lb}$ of $\mu_d$, its coefficient in $PD^{\mu_d}(F;u,v,t)$ (which is a polynomial  $c_{\lb}(u,v,t)$) can be recovered from the polynomial $c_{\lb}(u,v,-1)$. In other words, the monomials in $u,v$ coming from the various cohomology groups $H^j(F)$ should not undergo any simplication and their degree should tell from which cohomology group they come.

Consider first the trivial character $r_1$. Then $H^0(F)$ contributes to the coefficient $c_1(u,v,t)$ by 1 and $H^1(F)$ contributes by a multiple of the monomial $uvt$, since $H^1(F)_1=H^1(M)$ is still of pure  type $(1,1)$ in this more general setting. To see this, one may use the Gysin sequence
$$0=H^1(\PP^2\setminus \Sigma) \to H^1(M) \to H^0(C \setminus \Sigma)(-1) \to...$$
with $ \Sigma$ denoting the  set of singular points of the curve $C$.
The group $H^2(F)_1=H^2(M)$ has weights at least 2, since $M$ is smooth. On the other hand, the elements of weight 2 are those in the image of the morphism
$$i^*:H^2(\PP^2) \to H^2(M)$$
induced by the inclusion $i: M \to \PP^2$, since $\PP^2$ is a compactification of $M$. But this morphism is trivial,
since $H^2(\PP^2, \Q)$ is spanned by the first Chern class of the line bundle $\OO(d)$ and the restriction $\OO(d)|M$ is trivial. Indeed, it admits $Q$ as a section without zeroes. It follows that all the classes in $H^2(M)$ have in fact weights 3 and 4, and hence we can recover $c_1(u,v,t)$ from $c_1(u,v,-1)$.

Consider now a nontrivial character $r_{\lb}$, i.e. $\lb\ne 1$. Then $H^0(F)$ contributes to the coefficient $c_{\lb}(u,v,t)$ by 0 and $H^1(F)$ contributes by a linear form in  $ut,vt$, since $H^1(F)_{\ne 1}$ is still of pure of weight 1 in this more general setting, see  Theorem 1.5 in
 \cite{BDS} or Theorem 4.1 in \cite{DP}. The contribution of $H^2(F)$ to $c_{\lb}(u,v,t)$ is by monomials of the form
$u^ av^ bt^2$ with $a+b \geq 2$, since $F$ is a smooth variety. This implies again that we can recover $c_{\lb}(u,v,t)$ from $c_{\lb}(u,v,-1)$, which ends the proof of Proposition \ref{prop3}.

\begin{rk} 
\label{rk5} 
Note that the information contained in the polynomial $Sp(\A)$ is equivalent to the information contained in the specialization $HD^{\mu_d}(F;u,1)$, see \cite{DL13}. However, even if $Sp(\A)$ is known to be combinatorially determined by \cite{BS}, it is an {\it open question} if the same holds for the 
Hodge-Deligne polynomial $HD^{\mu_d}(F;u,v)$ of the Milnor fiber of a hyperplane arrangement. Moreover the specialization $HD^{\mu_d}(F;u,1)$ does not determine the Hodge-Deligne polynomial $HD^{\mu_d}(F;u,v)$, even in the case
of a line arrangement $\A$ having only double and triple points, since $Sp(\A)$ does not determine the Papadima-Suciu invariant $\be_3(\A)$ (which cancels out when we set $v=1$ in $HD^{\mu_d}(F;u,v)$). To have an explicit example, refer to Examples 5.4 and 5.5 in \cite{CS},
where the realizations of the configurations $(9_3)_1$ and $(9_3)_2$ are shown to have distinct $b_1(F)$'s. They have the same spectra by  Theorem \ref{jumpcoef}, having the same number of lines and triple points.

\end{rk}

\end{document}